\documentclass[12pt]{article}
\usepackage{theorem}
\usepackage{graphics}
\newtheorem{lemma}{Lemma}

\newtheorem{theorem}[lemma]{Theorem}

{\theorembodyfont{\upshape}}
{\theorembodyfont{\upshape}}
{\theorembodyfont{\upshape}}
{\theorembodyfont{\upshape}}
{\theorembodyfont{\upshape}}
{\theorembodyfont{\upshape}}
{\theorembodyfont{\upshape}}

\newcommand{\Z}{{\bf Z}}
\newcommand{\R}{{\bf R}}

\newcommand{\rme}{{\rm e}}
\newcommand{\rmd}{\,{\rm d}}

\newcommand{\cB}{{\cal B}}

\newcommand{\cG}{{\cal G}}
\newcommand{\cH}{{\cal H}}

\newcommand{\cL}{{\cal L}}
\newcommand{\cM}{{\cal M}}

\newcommand{\sig}{\sigma}
\newcommand{\alp}{\alpha}

\newcommand{\gam}{\gamma}
\newcommand{\lam}{\lambda}
\newcommand{\del}{\delta}
\newcommand{\eps}{\varepsilon}

\newcommand{\Dom}{{\rm Dom}}

\newcommand{\Spec}{{\rm Spec}}

\newcommand{\norm}{\Vert}
\renewcommand{\Re}{{\rm Re}\,}
\renewcommand{\Im}{{\rm Im}\,}

\newcommand{\Proof}{\underbar{Proof}{\hskip 0.1in}}

\newcommand{\lin}{{\rm lin}}
\newcommand{\Schrodinger}{Schr\"odinger }

\newcommand{\la}{{\langle}}
\newcommand{\ra}{{\rangle}}
\newcommand{\pr}{\prime}
\newcommand{\emp}[1]{{\it #1}}

\newcommand{\nopic}[1]{}

\newcommand{\ops}{one-parameter semigroup}

\newcommand{\be}{\begin{equation}}
\newcommand{\ee}{\end{equation}}
\newcommand{\choices}[1]{\left\{\begin{array}{ll} #1 \end{array}\right.}
\usepackage{hyperref}
\title{APPROXIMATING SEMIGROUPS\\  BY USING PSEUDOSPECTRA}
\author{E.B. Davies}
\date{18 March 2003}
\begin{document}
\maketitle

\section{Introduction}

Let $A$ be the generator of a \ops\ $T_t$ acting in a Hilbert
space $\cH$. We discuss the numerical computation of $T_tf$, or
equivalently the solution of the initial value problem
\be
f^\pr (t)=Af(t)\label{ivp}
\ee
given $f(0)=f$. This involves several problems. The first is that
the spectral mapping property may fail; that is one may have
\[
\Spec(T_t)\backslash \{0\}\not= \{\rme^{\lam t}:\lam\in\Spec(A)\}.
\]
In particular $\norm T_t\norm$ may grow faster than $\rme^{st}$ as
$t$ increases, where
\[
s=\sup\{\Re(\lam):\lam\in\Spec(A)\}.
\]
This problem is well known and has been studied from many points
of view, but it remains a difficulty, even if $A$ has discrete
spectrum, \cite{ebd4,davsim,nagel,ren,ren1,ren2,tre4}.

The second problem arises for differential operators, particularly
in several space dimensions, when the matrix approximations have
very high dimensions. Even if $A$ has a sparse matrix, $T_t$
generally has a full matrix, so storing the matrix entries is not
feasible. The obvious solution is to find a subspace of relatively
small dimension which contains most information of interest. One
might try to do this by taking the linear span of a finite number
of eigenvectors, those for which the real parts of the eigenvalues
are largest. Unfortunately experience shows that for many
non-self-adjoint operators $A$, the eigenvectors do not form a
basis; indeed the norms of the spectral projections often increase
exponentially fast according to natural orderings of the
eigenvalues. This forces one to be very cautious about assuming
that a spectral expansion of some given $f$ will yield useful
results.

This phenomenon is linked to the appearance of non-trivial
pseudospectra. When this happens the determination of more than a
small number of eigenvalues may become numerically impossible.
Even if theorems about the convergence of the eigenfunction
expansion of a general $f\in\cB$ subject to a resummation method
can be proved, they have limited use if most of the eigenvalues
and eigenvectors cannot be determined.

Many recent papers about pseudospectra have drawn attention to
possible instability problems which are not revealed by looking at
the spectrum alone,
\cite{tre1,tre2,website,wright,wright2,wright3}. Our goal in this
paper is more positive: we use pseudospectral methods to solve the
evolution equation above for highly non-self-adjoint operators.
The existence of a large number of approximate eigenvalues is
regarded as a resource rather than an embarrassment. We develop an
`approximate spectral expansion' which may have little to do with
the true eigenvalues and eigenvectors of the operator. In spite of
this our main result, Theorem~\ref{maintheorem}, may be used to
solve the evolution equation to a high degree of accuracy. In the
examples studied numerically we demonstrate that it is far more
accurate than the normal spectral expansion.

Our method is particularly useful if one wishes to solve the
initial value problem (\ref{ivp}) for a large number of different
choices of the initial data. The approximate eigenvalues and
eigenvectors only need to be produced once, and the computations
needed for each choice of the initial data are fairly easy.

The examples which we consider in this paper are
convection-diffusion operators. There are arguments in favour of
studying the associated semigroups $T_t=\rme^{At}$ in $L^1$ rather
than $L^2$. Diffusion is a probabilistic phenomenon, and the
conservation of probability is not easy to study in an $L^2$
context. It is shown in \cite{ebd4,davsim,murata2,zhang} that the
`same' semigroup may have different growth properties when studied
in $L^1$ or in $L^2$. Nevertheless we will focus on the $L^2$
theory, for the same reason as in classical Fourier theory: the
theorems are much simpler to state and apply.

\section{The Abstract Setting}

We start with several assumptions. The first is the choice of
numbers $M$, $\gam$ such that
\be
\norm T_t\norm \leq M\rme^{\gam t}\label{semigrowth}
\ee
for all $t\geq 0$. The second is the existence of a set $S$
equipped with a $\sig$-field of subsets and a finite measure $\rmd
s$. We assume that we are given measurable families of unit
vectors $u_s\in\cH$ and of complex numbers $\lam_s$ parametrized
by $s\in S$ and satisfying
\be
\inf\{\norm u_s -w\norm + \norm Aw-\lam_sw\norm
:w\in\Dom(A)\}<\eps. \label{generalized}
\ee
Throughout this paper $\eps$ is a given `acceptable' error
satisfying $0 <\eps < 1/2$. From a purely theoretical point of
view the assumption
\be
\norm Au_s-\lam_su_s\norm <\eps\label{epsbound}
\ee
for all $s\in S$ would be simpler. We prefer (\ref{generalized})
because it permits simpler expressions for the vectors $u_s$ in
applications. Clearly $\lam_s$ are approximate eigenvalues of $A$,
up to the error $\eps
>0$. The assumption (\ref{generalized}) implies that
\[
\lam_s\in\Spec_{2\eps}(A):=\{z:\norm (zI-A)^{-1}\norm\geq
(2\eps)^{-1}\}
\]
in the language of pseudospectral theory.

If $A$ is highly non-self-adjoint, the fact that $\lam_s$ are
approximate eigenvalues of $A$ does not imply that they are close
to the spectrum of $A$. This allows us to go far beyond what is
possible by means of conventional spectral methods. In numerical
applications we will take $S$ to be finite, but the above setting
allows a better understanding of the general theory.

We define a bounded, linear `pseudospectral' transform $\cG$ from
$L^1(S)$ to $\cH$ by
\[
\cG\phi=\int_S\phi(s)u_s\,\rmd s.
\]
We restrict $\cG$ to $L^2(S)$ and note that it is then bounded
with $\norm \cG\norm \leq |S|^{1/2}$, where $|S|$ is the measure
of $S$. The adjoint operator $\cG^\ast:\cH\to L^2(S)$ is given by
\[
(\cG^\ast f)(s)=\la f,u_s\ra
\]
and $B=\cG^\ast \cG:L^2(S)\to L^2(S)$ is given by
\[
(B\phi)(s)=\int_S b(s,t)\phi(t)\,\rmd t
\]
where
\[
b(s,t)=\la u_t,u_s\ra.
\]
Since $b$ is a bounded measurable function, $B$ is a
Hilbert-Schmidt operator on $L^2(S)$.

It is immediate from the definitions that $B\phi=0$ if and only if
$\cG\phi=0$. We assume throughout the paper that $B$ is
invertible, a matter which needs to be confirmed in any
application.

The following theorem describes how best to approximate $f\in\cH$
by expressions of the form $\cG \phi$ where $\phi\in L^2(S)$. We
will frequently refer to the algebraic sum $\cM=\cL+\cL^\perp$,
where $\cL$ is the range of $\cG$ in $\cH$. This is a dense linear
subspace of $\cH$. If $S$ is finite, as in all numerical
applications, then $\cL$ is closed and $\cM=\cH$.

\begin{theorem}\label{PN}
If $P$ is the orthogonal projection on $\cH$ with range
$\overline{\cL}$, then
\[
Pf=\cG B^{-1}\cG^\ast f
\]
for all $f\in\cM$. For such $f$ the quantity $\norm
\cG\phi-f\norm$, where $\phi\in L^2(S)$, is minimized by
$\phi=B^{-1}\cG^\ast f$. We also write $\phi=\cG\backslash f$, as
in Matlab.
\end{theorem}

\Proof If $f\in \cL^\perp$ then $\cG^\ast f=0$, so $\cG
B^{-1}\cG^\ast f=0$. If $f=\cG \phi$ then
\[
\cG B^{-1}\cG^\ast f= \cG (B^{-1}\cG^\ast\cG)\phi =\cG\phi =f.
\]
This proves the first statement. If $f=\cG\psi+g$ where $\psi\in
L^2(S)$ and $g\in\cL^\perp$ then
\[
\norm f-\cG\phi\norm^2=\norm \cG(\phi-\psi)\norm^2+\norm g\norm^2.
\]
This is clearly minimized for $\phi=\psi$ and, under our standing
hypothesis that $\cG$ is one-one, this is the unique minimum. We
also have
\[
B^{-1}\cG^\ast f=B^{-1}\cG^\ast (\cG\psi+g)=\psi.
\]

The above method of approximation should be contrasted with the
following alternative. Suppose that $A$ has a complete set of
eigenvectors $u_n$, $n=1,2,...$ and that $u_n^\ast$ are
corresponding eigenvectors of $A^\ast$, so that the two sets form
a biorthogonal system in the sense that $\la u_n,u_m^\ast\ra
=\del_{m,n}$. The standard spectral expansion with respect to this
system is
\be
f=\lim_{N\to\infty}Q_Nf  \label{basis}
\ee
where
\be%
Q_Nf= \sum_{n=1}^N \la f,u_n^\ast\ra u_n.\label{QN}
\ee%
If the identity (\ref{basis}) holds for all $f\in\cH$ one says
that $u_n$ form a basis in $\cH$. Unfortunately this is rarely
true for highly non-self-adjoint operators. Indeed $\norm
Q_N\norm$ frequently diverges at an exponential rate in
applications. One might modify the above formula by assuming
Cesaro or Abel summability, but convergence would still have to be
verified and is not always true.

One the other hand if the set $\{u_n\}_{n=1}^\infty$ is complete
we always have
\[
f=\lim_{N\to\infty} P_N f
\]
where $P_N$ is the orthogonal projection of $\cH$ onto
$\lin\{u_n:1\leq n\leq N\}$, and this is indeed the optimal
approximation sequence to $f$. Putting $S=\{1,...,N\}$ the theorem
above enables one to compute $P_N$. The main disadvantage of the
projections $P_N$ is that they do not commute with $A$.

Returning to the general context at the start of this section, we
use the operators defined above to solve the evolution equation
approximately. We start by obtaining a bound on the real parts of
approximate eigenvalues.

\begin{lemma}\label{appp}
If $\norm  u \norm =1$ and
\[
\norm u-w\norm +\norm Aw-\lam w\norm<\eps
\]
then $\Re(\lam)\leq \gam +2M\eps$.
\end{lemma}

\Proof We first observe that $0<\eps <1/2$ implies $1/2\leq \norm
w \norm \leq 3/2$. Putting $\mu=\Re(\lam)$, the identity
\[
\frac{\rmd }{\rmd s}\left\{T_{t-s}\rme^{\lam
s}w\right\}=T_{t-s}\rme^{\lam s}(\lam w-Aw)
\]
implies that
\begin{eqnarray*}
\norm T_tw-\rme^{\lam t}w\norm &=& \norm \int_0^t
T_{t-s}\rme^{\lam s}(\lam w-Aw)\,\rmd s\norm\\
&\leq & \int_0^tM\rme^{\gam(t-s)+\mu s}\eps\,\rmd s\\
&=&M\eps\frac{  \rme^{\mu t}-\rme^{\gam t} }{ \mu-\gam}
\end{eqnarray*}
It follows that
\[
\rme^{\mu t}/2\leq 2M\rme^{\gam t}+ M\eps\frac{  \rme^{\mu
t}-\rme^{\gam t} }{ \mu-\gam}.
\]
If $\mu >\gam$ then letting $t\to +\infty$, we deduce that
\[
1\leq \frac{2 M\eps }{\mu- \gam }
\]
which is equivalent to the statement of the lemma.

It follows immediately from the lemma that if $\Re(\lam)>\gam$ and
we put $\tilde{\lam}=\gam +i\,\Im(\lam)$ then
$|\lam-\tilde{\lam}|\leq 2M\eps$. Therefore
\[
\norm u-w\norm + \norm Aw-\tilde{\lam} w\norm<\eps(3M+1).
\]
In the rest of the paper we assume that these changes in the
approximate eigenvalues have been made, so that $\Re(\lam_s)\leq
\gam$ for all $s\in S$, and that $\eps$ has been increased
correspondingly.

The main theorem of this paper is best formulated in terms of
certain approximating semigroups $R_t$.

\begin{theorem}\label{localapprox}
Let $\rme^{\lam t}$ be the multiplication operator on $L^2(S)$
defined by
\[
(\rme^{\lam t}\phi)(s)=\rme^{\lam_s t}\phi(s).
\]
and define $R_t$ on $\cM$ for $t\geq 0$ by
\be
R_t=\cG\rme^{\lam t} B^{-1}\cG^\ast\label{approxsemigroup}
\ee
Then $R_0=P$, $R_t(\cL^\perp)=0$ and $R_t(\cL)\subseteq \cL$ for
all $t\geq 0$. We also have
\[
R_tR_uf=R_{t+u}f
\]
for all $t,u\geq 0$ and $f\in \cM$.
\end{theorem}

\Proof If $f\in \cL^\perp$ then $\cG^\ast f=0$ so $R_tf=0$. If
$f=\cG\phi$ where $\phi\in L^2(S)$ then
\[
R_tf=\cG\rme^{\lam t} B^{-1}\cG^\ast\cG\phi=\cG(\rme^{\lam t}
\phi)\in\cL.
\]
Finally, if $f=\cG\phi$ then
\[
R_tR_uf=R_t\cG(\rme^{\lam u}\phi)=\cG(\rme^{\lam t}\rme^{\lam
u}\phi)=R_{t+u}f.
\]

\begin{theorem} Suppose that $S$ is finite and $\cH=L^2(X,\rmd
x)$, and rewrite $u_s(x)=u(x,s)$. Then
\[
(R_tf)(x)=\int_X K_t(x,y)f(y) \, \rmd y
\]
for all $f\in \cH$, where
\[
K_t(x,y)=\sum_{r,s} u(x,s)\rme^{\lam_s
t}(B^{-1})_{s,r}\overline{u(y,r)}
\]
\end{theorem}

\Proof Since $S$ is finite, $\cL$ is a finite-dimensional subspace
of $\cH$, and $\cL +\cL^\perp=\cH$. We deduce that $R_t$ has
domain $\cH$. The formulae of the theorem are the result of
rewriting (\ref{approxsemigroup}) in integral operator form.

One might conjecture that the integral kernel of $R_t$ is
uniformly close to that of $T_t$ under suitable conditions, but we
do not have any such results.

The following is our main theorem. It is only numerically
efficient if $\eps>0$ and $\del=\norm f-Pf\norm$ are both small.
We discuss this further in the next section.

\begin{theorem}\label{maintheorem}
If $f\in\cM$ then
\be
\norm T_t f-R_tf\norm < \norm f-Pf\norm M \rme^{\gam t}+ \eps
(1+M+Mt)\norm \cG\backslash f\norm _1\rme^{\gam
t}\label{mainestimate}
\ee
for all $t\geq 0$.
\end{theorem}

\Proof If we put $\phi=\cG\backslash f=B^{-1}\cG^\ast f$ then the
estimate can be rewritten in the form
\be
\norm T_t f-\cG\phi_t\norm < \norm f-\cG\phi\norm M \rme^{\gam t}+
\eps (1+M+Mt)\norm \phi\norm _1\rme^{\gam t}\label{modified}
\ee
where $\phi_t=\rme^{\lam t}\phi$. We follow the argument of
Lemma~\ref{appp} up to
\begin{eqnarray*}
\norm T_tw-\rme^{\lam t}w\norm
&\leq & \int_0^tM\rme^{\gam(t-s)+\mu s}\eps\,\rmd s\\
&\leq& \int_0^tM\rme^{\gam t}\eps\,\rmd s\\
&=&\eps Mt\rme^{\gam t}.
\end{eqnarray*}
Hence
\begin{eqnarray*}
\norm T_tu-\rme^{\lam t}u\norm &\leq &\norm \rme^{\lam
t}(w-u)\norm +
\norm T_t(u-w)\norm + \norm T_tw-\rme^{\lam t}w\norm \\
&=&\eps (1+M+Mt)\rme^{\gam t}.
\end{eqnarray*}

Applying this to each $u_s$ in the expansion
\[
\cG\phi=\int_S \phi(s)u_s \, \rmd s
\]
yields
\[
\norm T_t\cG \phi-\cG\phi_t\norm \leq \eps (1+M+Mt)\norm
\phi\norm_1\rme^{\gam t}.
\]
The theorem follows by combining this with the bound
\[
\norm T_t f-T_t\cG\phi\norm \leq M\rme^{\gam t} \norm
f-\cG\phi\norm.
\]

The above theorem is only useful as long as the right hand side of
(\ref{mainestimate}) is much smaller than $\norm R_tf\norm$. Since
\[
\norm R_tf\norm=\norm \cG\phi_t\norm\leq\norm \phi\norm_1\rme^{\mu
t}
\]
where
\[
\mu=\sup\{\Re(\lam_s):s\in S\},
\]
the estimates are only useful for a short time if $\mu \ll \gam$.
The point here is that $\mu$ may be substantially larger than
$\sup\{\Re(z):z\in \Spec(A)\}$, so pseudospectral methods may be
correspondingly more accurate than spectral methods.

\section{Numerical implementation}

In numerical applications we take $S$ to be a finite set, possibly
containing fewer than a hundred points. This implies that
$\cM=\cH$. The main task is the choice of the vectors $u_s\in\cH$.
Once this has been done, there are three possible methods of
computing $\phi=B^{-1}\cG^\ast f$ given $f\in\cH$. The vectors
$u_s$ determine $\cG^\ast$, and also the operator $B$ via the
kernel $b(s,t)$. One might compute $B^{-1}$ and then apply the
above formula to obtain $\phi$. Since the operator $B^{-1}$ is
highly singular it is better to evaluate $B^{-1}\psi$ for
$\psi=\cG^\ast f$ without computing $B^{-1}$; Matlab uses the
command $B\backslash \psi$ for this purpose. One may finally avoid
any reference to $\cG^\ast$ or $B$, by using Matlab to compute
$\phi=\cG\backslash f$ directly. Since $\cG$ is a rectangular
matrix, Matlab actually finds the `solution' with least squares
error. We tried all three methods, and found, as expected, that
the third is by far the most accurate. Once $\phi$ has been
determined we do not use Theorem~\ref{maintheorem} as stated, but
the reformulation in (\ref{modified}).

The choice of a suitably small $\eps>0$ is made before starting
the computation. On the other hand the verification that $\norm
f-\cG\phi\norm$ is small is done on a posterior basis. Since
$\phi$ and $\cG$ have to be computed in any case, this poses no
problems.

There are two obvious ways of choosing the unit vectors $u_s$. One
may use one's physical intuition, as in the examples of this
paper, to select certain vectors, and then show that they satisfy
the fundamental inequality (\ref{generalized}) for a suitably
small $\eps >0$. This method has been used successfully in the
semiclassical context, \cite{semiclass1,semiclass2,DSZ,zwor}. The
second method depends upon numerical, pseudospectral calculations,
and will be described in more detail in a later publication. The
first stage is the replacement of the differential operator $A$ by
a sparse matrix approximation, possibly in a space of very high
dimension. This may involve finite element methods or wavelets,
and is not the focus of this article. There is now a
well-developed technology for calculating pseudospectra, and it
may be applied to very large sparse matrices. Given $\eps
>0$, we next have to choose a finite set of numbers from the set
$\Spec_\eps(A)$. If $t>0$ is known, there is no need to consider
points in $\lam\in\Spec_\eps(A)$ for which $\rme^{\lam t}$ is
extremely small, because the contributions of the corresponding
terms of $\cG\phi_t$ will be negligible. This applies in
particular to any eigenvalues of $A$ whose real parts are much
less than $\gam$. For each $\lam_s$ we finally choose a unit
vector $u_s$ for which $\norm Au_s-\lam u_s\norm<\eps$.

In some cases it might be advisable to chose several vectors $u_s$
corresponding to each $\lam_s$, providing each vector with a
different label $s$. The choice depends upon how many eigenvalues
of order $\eps^2$ the operator
\[
D_s=(\lam_sI-A)^\ast(\lam_sI-A)
\]
possesses. For rotationally invariant problems in dimension two,
for example, one would treat each angular momentum sector
independently, and include that parameter in the labelling of $S$.
The $\eps$-pseudospectra for different sectors may well overlap.

\section{A Pure Convection Operator}

The theory above has applications to convection-diffusion
operators, but the simplest example is given by the pure
convection operator
\[
(Af)(x)=f^\pr(x)
\]
acting in $L^2(0,a)$ subject to the boundary condition $f(a)=0$.
This is the generator of the \ops\ $T_t$ given by
\[
(T_tf)(x)=\choices{f(x+t)&\mbox{ if $x+t<a$}\\
0&\mbox{ otherwise.} }
\]
Since $T_t=0$ for all $t\geq a$, $\gam$ can take any value in the
estimate (\ref{semigrowth}). Nevertheless, since we are interested
primarily in the case of large $a$, we take $\gam=0$ and $M=1$.
The fact that $\Spec(A)=\emptyset$ implies that one cannot hope to
use spectral expansions to evaluate $T_t$, but pseudospectral
expansions are still possible. Since this example is exactly
soluble, we only analyze it by our method in order to understand
how well the method works. We will see in the next section that
the pseudospectral expansion of this operator is an asymptotic
form of the corresponding expansion for a simple
convection-diffusion operator.

The following constructions depend upon the choice of positive
constants $c$ and $\alp$. In many cases an appropriate value of
$c$ may be found in the range $5\leq c\leq 10$. The value $c=0$
leads to a Fourier series expansion, which is not appropriate for
this problem. One could put $\alp=1$, but for asymptotic theorems
it might be more appropriate to make it proportional to $a$ and/or
inversely proportional to $c$. Let $v:[0,1]\to [0,1]$ be the
function
\[
v(x)=\choices{ 1&\mbox{if $0\leq x\leq a-\alp$}\\
(a-x)/\alp&\mbox{if $a-\alp\leq x\leq a$.} }
\]
(Many other choices would be equally suitable, for example
$v(x)=1-\rme^{(x-a)/\alp}$.) Given $s\in\Z$ we define $u_s\in
L^2(0,a)$ by
\[
u_s(x)=k\rme^{-cx/a+2\pi isx/a}
\]
where
\[
k^{-2}=\frac{a}{2c}\left\{ 1-\rme^{-2c}\right\}.
\]
This choice implies the identity $\norm u_s\norm =1$. We also see
that $ak^2/2c\to 1$ at an exponential rate as $c$ increases. If we
define $w_s\in\Dom(A)$ by
\[
w_s(x)=u_s(x)v(x)
\]
then
\begin{eqnarray*}
\norm u_s-w_s\norm^2&=& k^2\int_0^a \rme^{-2cx/a} (1-v(x))^2\,\rmd x\\
&\leq & \frac{ak^2}{2c}\rme^{-2c(1-\alp/a)}.
\end{eqnarray*}
If we put $\lam_s=-c/a+2\pi is/a$ then
\begin{eqnarray*}
\norm Aw_s-\lam_sw_s\norm^2 &=&k^2\int_0^a \rme^{-2cx/a} v^\pr(x)^2\,\rmd x\\
&\leq & \frac{ak^2}{2c\alp^2}\rme^{-2c(1-\alp/a)}.
\end{eqnarray*}
This indicates that the bound (\ref{generalized}) holds with
$\eps=O(\rme^{-c(1-\alp/a)})$ as $c\to\infty$.

Having chosen a sufficiently large $N >0$, we then put
\be
S=\{s\in\Z :-N \leq s\leq N\}.\label{intervalS}
\ee

The integral kernel of $B=\cG^\ast\cG$ is
\begin{eqnarray*}
b(s,t)&=&\la u_t,u_s\ra\\
&=& k^2\int_0^a\rme^{-2cx/a +2\pi i(t-s)x/a} \,\rmd x\\
&=&k^2\frac{1-\rme^{-2c}}{2c/a-2\pi i(t-s)/a} \\
 &\sim &\frac{1}{1-\pi i(t-s)/c}
\end{eqnarray*}
if $c$ is sufficiently large.

If $f\in L^2(0,a)$ and $g=\cG^\ast f$ then
\[
g(s)= k\int_0^{a} f(x)\rme^{-cx/a-2\pi isx/a}\,\rmd x
\]
and
\be
(Pf)(x)= k\rme^{-cx/a}\sum_{s=-N}^N (B^{-1}g)(s)\rme^{2\pi isx/a}
.\label{newfourier}
\ee
If $Pf$ is approximately equal to $f$ then we have shown that
\be
(T_tf)(x)\sim k\rme^{-cx/a}\sum_{s=-N}^N (B^{-1}g)(s)\rme^{2\pi
isx/a+(-c/a+2\pi i sx/a) t}\label{newsemi}
\ee
for $t>0$. In numerical implementations one actually uses the
equivalent formula
\be
(T_tf)(x)\sim k\rme^{-cx/a}\sum_{s=-N}^N \phi(s)\rme^{2\pi
isx/a+(-c/a+2\pi i sx/a) t}\label{compsemi}
\ee%
where $\phi=\cG\backslash f$, in the notation of Matlab.

Let us compare this with what one gets by using ordinary Fourier
series, by making the choices $c=0$ and $N=\infty$ in the above
formulae. We then have $k=a^{-1/2}$ and
\[
u_s(x)=a^{-1/2}\rme^{2\pi i s x/a}
\]
for all $s\in\Z$. We also have
\[
g(s)=a^{-1/2}\int_0^a f(x)\rme^{-2\pi i s x/a}\,\rmd x
\]
so $g$ is the sequence of Fourier coefficients of $f$, assuming
periodic boundary conditions. Since
\[
b(s,t)=\choices{1&\mbox{if $s=t$}\\
0&\mbox{otherwise,}}
\]
$B$ is the identity operator on $l^2(\Z)$, and (\ref{newfourier})
is replaced by
\[
f(x)=a^{-1/2}\sum_{s=-\infty}^\infty g(s)\rme^{2\pi i sx/a},
\]
while (\ref{newsemi}) is replaced by
\begin{eqnarray*}
(\tilde{T_t}f)(x)&=&a^{-1/2}\sum_{s=-\infty}^\infty g(s)\rme^{2\pi
i sx/a+2\pi i st/a}\\
&=& f(x+t),
\end{eqnarray*}
subject to \emp{periodic boundary conditions} on $[0,a]$. The use
of Fourier series therefore solves a different problem from that
in which we are interested.

We implemented the above ideas numerically for two choices of the
initial function. We put $a=20$ and divided each unit interval
into $50$ equally spaced points, so that functions on $[0,a]$ are
approximated by sequences with $1000$ terms. We chose the initial
function to be
\[
f(x)=2\rme^{-10(x-5)^2}-\rme^{(x-5)^2/10}.
\]
We defined $f_t$ to be the right-hand side of (\ref{newsemi}) and
computed
\[
p=\norm f-f_0\norm_\infty\hspace{1 cm} q=\norm
T_tf-f_t\norm_\infty
\]
for various values of $c,\, N$, putting $t=5$. (Similar results
are obtained using the $L^2$ norm.)  The results are presented in
Table 1.
\begin{table}[h]
\[
\begin{array}{cccc}
c&N&p&q\\
5 & 30  & 0.056  &  0.056 \\
10& 30  & 0.049  & 0.049  \\
5 &  40 & 0.0075  &  0.0075 \\
10&  40 &  0.0063 &  0.0063 \\
3 & 50  &  0.0040 & 0.050  \\
5 & 50  &  0.00063 & 0.0067  \\
10& 50  &  0.00051 & 0.00051  \\
\end{array}
\]
\begin{center}
Table 1
\end{center}
\end{table}
Our conclusion from the data is that the errors depend more upon
the number of terms $2N+1$ in the expansion than upon the value of
$c$. However, in the best case, $N=50$, we see that $c$ needs to
be substantially bigger than $5$ for accurate results.

We also considered the initial function $f=1$, for which $T_tf$ is
the characteristic function of $(0,a-t)$. We made the same choices
$a=20$, $N=50$, $c=10$ and $t=5$ as above. This case is highly
singular, since neither $f$ nor $T_tf$ are close to being in the
domain of $A$. Although the computed values of $f_t$ are close to
$1$ for $x<15$ and close to $0$ for $x>15$, there is a Gibbs-type
phenomenon near $x=15$, the maximum value of $f_t$ being about
$1.21$. As expected, the maximum is unchanged for $N=100$.

The example of this section may be described in terms of a global
approximating semigroup, to be contrasted with the local
approximating semigroups of Theorem~\ref{localapprox}. We
introduce the operator
\[
(A_cf)(x)=f^\pr(x)
\]
acting in $L^2(0,a)$ subject to the boundary conditions
$f(a)=\rme^{-c}f(0)$. If $c\geq 0$ this is the generator of a
\ops\ $T_{c,t}$ acting on $L^2(0,a)$. One sees immediately that
$u_s,\, \lam_s$ are the eigenvectors and eigenvalues respectively
of $A_c$. Section 2 provides estimates of how closely solutions of
$f^\pr(t)=Af(t)$ are approximated by solutions of
$f^\pr(t)=A_cf(t)$ which involve only a finite number of
eigenvectors of $A_c$. However, the right-hand side of
(\ref{newsemi}) is not simply a spectral expansion of $T_{c,t}$.
The similarity between $T_t$ and $T_{c,t}$ explains why we should
expect steadily better approximations as $c$ increases, provided
the computations remain feasible.

\section{A Convection-Diffusion Operator}

The difference between the $L^1$ and $L^2$ behaviour of semigroups
is well illustrated by the pure convection operator
\[
(Af)(x)=-2xf^\pr (x)
\]
which generates the semigroup
\[
(T_tf)(x)=f(\rme^{-2t}x).
\]
This is a positivity preserving contraction semigroup on
$C_0(\R)$, but on $L^2(\R)$ we have $\norm T_t\norm=\rme^t$ for
all $t\geq 0$. The semigroup has the same behaviour when acting on
$L^2(-a,a)$, and if $b >0$ is large enough one would expect
similar behaviour for the convection-diffusion operator
\[
(Af)(x)=b^{-1} f^{\pr\pr}(x)-2xf^\pr(x)
\]
acting in $L^2(-a,a)$ subject to Dirichlet boundary conditions at
$\pm a$.

We consider the somewhat simpler operator
\[
(Af)(x)=b^{-1} f^{\pr\pr}(x)+f^\pr(x)
\]
acting in $L^2(0,a)$ in more detail. The first term produces a
diffusion effect while the second cause a drift to the left at
speed $1$. If we impose Dirichlet boundary conditions then
$T_t=\rme^{At}$ is a positivity preserving contraction semigroup
on $L^p(0,a)$ for all $1\leq p<\infty$.  If $b$ is large then the
norm of $T_t$ remains close to $1$ for $t$ up to about $a$ and
then decreases rapidly towards $0$. We put $M=1$ and $\gam=0$ in
our theorems.

The following results are well-known, \cite{RT,ebd3}. The
eigenvectors and eigenvalues of $A$ are given by
$e_n(x)=k_n\rme^{-bx/2}\sin(\pi n x/a)$ and $\lam_n=-b/4-\pi^2
n^2/ba^2$ respectively for $n=1,2,...$. We have $\norm e_n \norm
=1$ for all $n$ if
\begin{eqnarray*}
k_n^{-2}&=& \frac{1}{4}\int_0^a \rme^{-bx}\left| \rme^{\pi i
nx/a}-\rme^{-\pi i nx/a}\right|^2\, \rmd x\\
&=&\frac{2\pi^2 n^2(1-\rme^{-ba})}{b(b^2a^2+4\pi^2n^2)}.
\end{eqnarray*}
We see that
\be%
k_n^2\sim \frac{b^3a^2}{2\pi^2n^2}\label{kasympt}
\ee%
as $b\to\infty$ for fixed $n,\, a$. The spectrum of $A$ is
asymptotically empty as $b\to\infty$ for fixed $a$. It converges
to $(-\infty, -b/4]$ as $a\to \infty$, but this is not the
spectrum of $A$ considered either in $L^2(0,\infty)$ or in
$L^2(\R)$. The normalized eigenvectors of $A^\ast$ are $e_n^\ast=
k_n\rme^{b(x-a)/2}\sin(\pi n x/a)$.

\begin{lemma}
The two sets of eigenvectors $\{ e_n\}$ and $\{e_m^\ast\}$ satisfy
$\la e_n,e_m^\ast\ra =0$ if $m\not= n$. The corresponding spectral
projections $P_n$ of $A$ satisfy
\[
\norm P_n \norm \sim\frac{2\pi^2n^2}{b^3a^3}\rme^{ba/2}
\]
as $b\to\infty$ for each $n,\, a$.
\end{lemma}

\Proof The first statement can be verified directly, but it is a
consequence of the fact that the two sets are eigenvectors of $A$
and $A^\ast$ respectively. A direct calculation shows that
\be%
\la e_n, e_n^\ast\ra =\frac{k_n^2 a}{2\rme^{ba/2}}.\label{inner}
\ee%
The second statement now follows by substituting (\ref{kasympt})
and (\ref{inner}) into
\[
\norm P_n \norm = |\la e_n, e_n^\ast\ra|^{-1}.
\]

All of the above facts suggest that one should not use spectral
expansions for large $b$.

In order to test this we computed $P_Nf$
 as defined by Theorem~\ref{PN} with
$\cL=\lin\{e_1,...,e_N\}$ and $Q_Nf$ as defined by (\ref{QN}). We
chose $a=20$ and discretized using $10$ points per unit interval,
so that $[0,a]$ was replaced by a set of $201$ points, including
the endpoints. We took the function $f$ to be
\[
f(x)=\rme^{-(x-a/2)^2}.
\]
Table 2 shows the sizes of  $ p=\norm f-P_Nf\norm$ and $q=\norm
f-Q_Nf\norm$ for a range of choices of $N$ when $b=2.5$ and when
$b=5.0$. We see that both methods have comparable accuracy for
$N=100$. However, the method using $P_N$ attains this accuracy far
more rapidly as $N$ increases than the pure spectral method using
$Q_N$. As $b$ increases the convergence of both methods
deteriorates, and for $b=7.5$ neither method gives useful results
for any value of $N$ up to $100$.

\begin{table}[h]
\[
\begin{array}{ccc}
b=2.5&\hspace*{10mm}&b=5.0\\

\begin{array}{ccc}
    N   &    p     &    q\\
10 &  3.9\times 10^{-1}    & 1.8\times 10^{3}   \\

   20  & 4.6\times 10^{-2}  & 1.3 \times 10^{3} \\

   30 &  1.7\times 10^{-3}  &  5.7\times 10^{1} \\

   40 &  1.8\times 10^{-5}  &  1.4\times 10^{-1} \\

   50 &  5.3\times 10^{-8}  &  2.0\times 10^{-3} \\

   60 &  4.1 \times 10^{-11}  &  2.9\times 10^{-6} \\

   70  & 1.7 \times 10^{-11}  &  4.4\times 10^{-10} \\

   80 &  9.0 \times 10^{-12}  &  1.4\times 10^{-10} \\

   90 &  1.2 \times 10^{-11}  &  2.0\times 10^{-10} \\

   100 &    9.3\times 10^{-12}  & 3.4 \times 10^{-10} \\
   \end{array}
&&
\begin{array}{ccc}
    N   &    p     &    q\\
 10 &  7.3\times 10^{-1}    & 3.3\times 10^{9}   \\

   20  & 1.1\times 10^{-1}  &  1.4\times 10^{7} \\

   30 &  4.8\times 10^{-3}  &  3.3\times 10^{7} \\

   40 &  1.4\times 10^{-4}  &  2.8\times 10^{5} \\

   50 &  1.9\times 10^{-4}  &  8.3\times 10^{2} \\

   60 &  3.4 \times 10^{-5}  & 1.7 \times 10^{0} \\

   70  &  4.8\times 10^{-4}  &  6.0\times 10^{-5} \\

   80 &  2.4 \times 10^{-5}  &  5.9\times 10^{-5} \\

   90 &  1.4 \times 10^{-5}  &  9.5\times 10^{-5} \\

   100 &   1.0 \times 10^{-4}  & 2.1 \times 10^{-4} \\

   \end{array}
   \end{array}
\]
\begin{center}
Table 2
\end{center}
\end{table}


Our goal in the remainder of this section is to demonstrate that
pseudospectral expansions are useful for much larger values of
$b$. For any choice of $b$ the pseudospectra behave in an
interesting way as $a$ increases. For every $z$ inside the
parabola $\sig\in\R \to -b^{-1}\sig^2+i\sig$ one has
\[
\lim_{a\to\infty}\norm (zI-A)^{-1}\norm =+\infty
\]
and one can construct approximate eigenfunctions for all such $z$
by the following method. Given $\del$ satisfying $0<\del<1/2$ and
$\sig\in\R$, we put
\be%
u_\sig(x)=k\left(
\rme^{(-b/2+b\del+i\sig)x}-\rme^{(-b/2-b\del-i\sig)x}\right)
\label{usig}
\ee%
where $k=k(b,\del,\sig,a)$ is given by
\[
k^{-2}= \int_0^a\left|\rme^{(-b/2+b\del+i\sig)x}
-\rme^{(-b/2-b\del-i\sig)x}\right|^2\,\rmd x.
\]
Clearly $\norm u_\sig\norm =1$. We make $\del$ depend upon $a$
according to the formula
\[
\del=1/2-c/(ab)
\]
where $0<c<ab/2$. (As before one might choose $c$ in the range
$5\leq c \leq 10$.) This choice of $\del$ ensures that $k^{-2}\sim
a(1-\rme^{-2c})/2c$, $|u_\sig(a)|\sim k\rme^{-c}$, $u_\sig(0)=0$
and $u_\sig^\pr (0)\sim k(b+2i\sig)$ as $a\to \infty$. We also put
\[
w_\sig(x)=u_\sig(x)v(x)
\]
where $v(x)=1-\rme^{(x-a)/\alp}$, and $\alp$ is a constant such as
$\alp=1$. Finally we put
\begin{eqnarray*}
\mu_\sig&=&b^{-1}(-b/2+b\del+i\sig)^2+(-b/2+b\del+i\sig)\\
&=&(b\del^2-b/4)-b^{-1}\sig^2+2i\del \sig\\
&=& -b^{-1}\sig^2+i\sig-c/a+c^2/(a^2b) -2i\sig c/(ab)\\
&\to & -b^{-1}\sig^2+i\sig
\end{eqnarray*}
as $a\to\infty$.

There is no reason to expect that taking a large value of $b$
should cause problems. Indeed, as $b\to\infty$, the functions
$u_\sig$ defined by (\ref{usig}) converge to the corresponding
functions $u_s$ defined for the pure convection operator of
Section 4. In both cases the size of the constant $c$ controls the
degree of accuracy of the fundamental estimate
(\ref{generalized}). As we have seen before, this has to be
weighed against the increased difficulty of performing the
computations for large $c$.

\begin{theorem}
Under the above conditions there exists a constant
$K_{\alp,b,\sig}$ such that
\[
\norm u_\sig -w_\sig\norm + \norm Aw_\sig-\mu_\sig w_\sig\norm
\leq K_{\alp,b,\sig}
a^{-1/2}\left\{\frac{2c}{\rme^{2c}-1}\right\}^{1/2}
\]
for large enough $a>0$.
\end{theorem}

\Proof We have
\begin{eqnarray*}
\norm u_\sig-w_\sig\norm^2&=&k^2\int_0^a \left|
\rme^{(-b/2+b\del+i\sig)x}-
\rme^{(-b/2-b\del-i\sig)x}\right|^2\rme^{2(x-a)/\alp}\,\rmd x\\
&\leq &4k^2\int_0^a \rme^{(-b+2b\del)x+2(x-a)/\alp}\,\rmd x\\
&=&\frac{2k^2}{1/\alp-c/a} \left( \rme^{-2c}-\rme^{-2a/\alp}\right)\\
&\leq & 3k^2\alp\rme^{-2c}
\end{eqnarray*}
for large enough $a>0$.

Since $b^{-1}u_\sig^{\pr\pr}+u_\sig^\pr=\mu_\sig u_\sig$ and
$w_\sig\in\Dom(A)$, we have
\[
Aw_\sig-\mu_\sig w_\sig= 2b^{-1}u_\sig^\pr v^\pr+b^{-1}u_\sig
v^{\pr\pr}+u_\sig v^\pr.
\]
Therefore
\[
\norm Aw_\sig-\mu_\sig w_\sig\norm \leq 2b^{-1}\norm u_\sig^\pr
v^\pr\norm +b^{-1}\norm u_\sig v^{\pr\pr}\norm +\norm u_\sig
v^\pr\norm .
\]
Each of the terms on the right-hand side is estimated in the same
way as above. For example
\begin{eqnarray*}
\norm u_\sig v^\pr\norm^2&=&k^2\alp^{-2}\int_0^a \left|
\rme^{(-b/2+b\del+i\sig)x}-
\rme^{(-b/2-b\del-i\sig)x}\right|^2\rme^{2(x-a)/\alp}\,\rmd x\\
&\leq &\frac{2k^2}{\alp^2(1/\alp-c/a)} \left( \rme^{-2c}-\rme^{-2a/\alp}\right)\\
&\leq & 3k^2\alp^{-1}\rme^{-2c}
\end{eqnarray*}
for large enough $a>0$. Combining all these estimates yields the
statement of the theorem.

For general values of $\sig\in\R$ the functions $u_\sig$ do not
satisfy any set of linear boundary conditions. However, if we put
$\sig=2\pi s/a$ where $s\in\Z$ then there exist non-zero constants
$c_i$ such that $u_\sig(0)=0$, $u_\sig^\pr(0)=c_1+c_2 \sig$,
$u_\sig (a)=c_3$ and $u_\sig^\pr(a)=c_4+c_5\sig$. Therefore the
functions $u_\sig$ all satisfy boundary conditions of the form
$u(0)=0$ and
\[
c_5u^\pr(0)-c_2 u^\pr(a)=c_6u(a).
\]

We tested the above ideas numerically. We re-parametrized by means
of the substitution $\sig=2\pi s/a$ where $s\in\Z$ and $-N\leq
s\leq N$. We put $a=20$, each unit interval in $[0,a]$ being
represented by $10$ equally spaced points. We put $\alp=1$, $b=20$
and $c=5$. We took the same function $f$ as before, that is
\[
f(x)=\rme^{-(x-a/2)^2}.
\]
Table 3 shows the values of $p=\norm f-P_Nf\norm$ for various
values of $N$. The dimension of the subspace $\cL$ is $2N+1$.

\begin{table}[h]  
\[
\begin{array}{ccccc}
2N+1 &  p\\
11& 3.3\times 10^{-1}  \\
21&  3.4\times 10^{-2}  \\
31&  1.1\times 10^{-3}  \\
41&  1.1\times 10^{-5}  \\
51&  3.4\times 10^{-8}  \\
61&   2.9\times 10^{-11} \\
71&  1.3\times 10^{-14}\\
\end{array}
\]
\begin{center}
Table 3
\end{center}
\end{table}
The superiority of this method of expansion over both of the
previous ones is immediately clear. Further computations show that
the pseudo-spectral method works just as well for all values of
$b$ from $5$ to $100$ (and probably beyond that).

We finally computed the approximation $f_t=\cG\phi_t$ to $T_tf$
given by the formula (\ref{mainestimate}) of
Theorem~\ref{maintheorem}. We chose the parameters and initial
value of $f$ as above but put $N=15$; the choice $N=30$ gave the
same results up to the accuracy displayed. We discovered, as
expected, that $f_t$ is approximately non-negative; in fact
\[
-3\times 10^{-4}\leq\min\{f_t(x):0\leq x\leq a\}\leq 0
\]
for all $0\leq t\leq 16$, at which point we stopped the
computation. The shape of $f_t$ remained approximately gaussian as
$t$ increased, with the centre moving to the left and the width
slowly increasing. The maximum of $f_t$ decreases slowly up to
$t\sim 10$, when the centre of the peak approaches the origin,
after which it decreases rapidly. The graphs of $f$, $f_4$, $f_8$,
$f_{12}$ are plotted in Figure 1.

\begin{figure}[h]  
\begin{picture}(200,350)(0,0)
\includegraphics{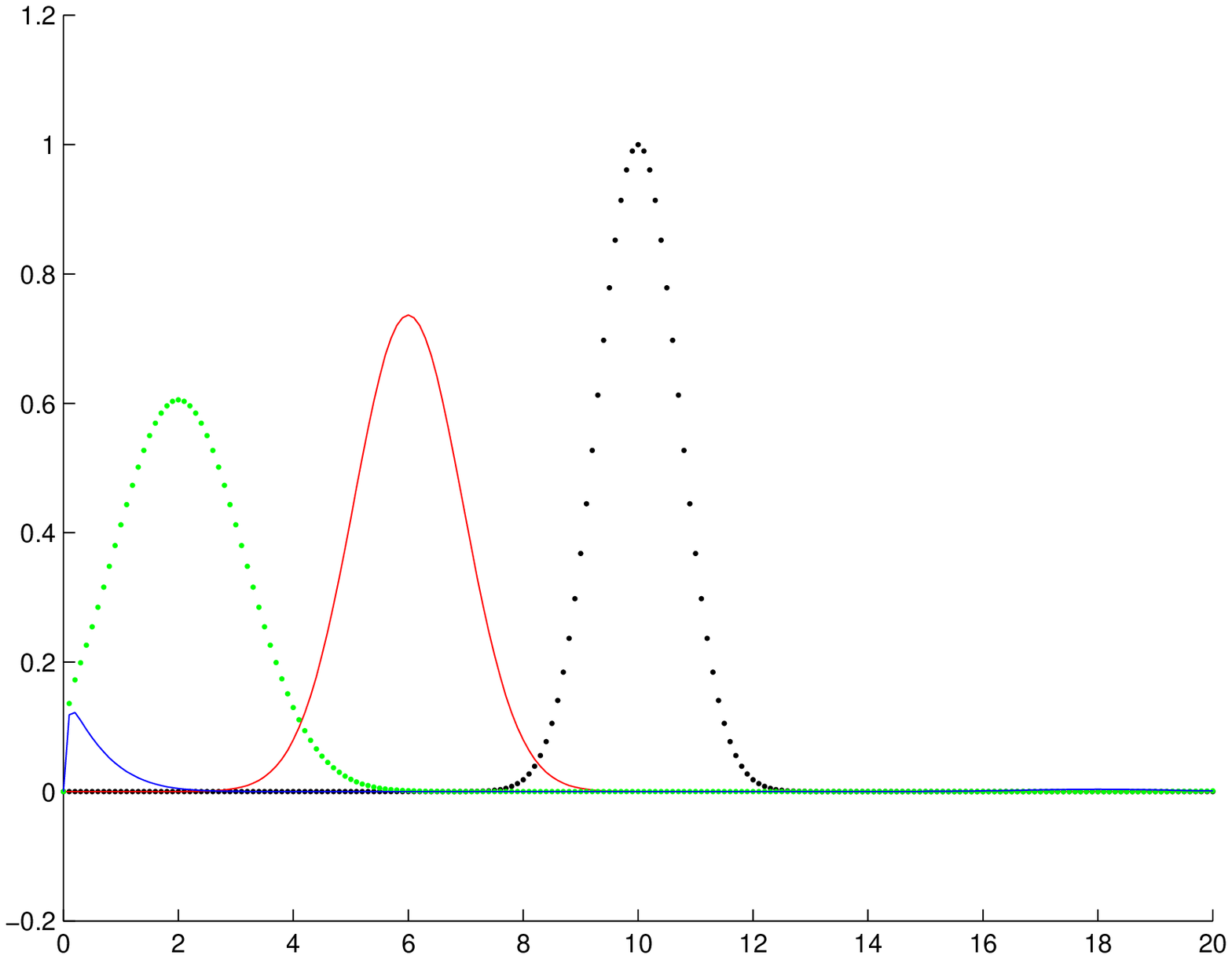}
\end{picture}
\begin{center}
Figure 1. Graphs of $f$, $f_8$ (dotted) and $f_4$, $f_{12}$
(solid)
\end{center}
\end{figure}
\clearpage

The detailed behaviour of the maximum $m$ is presented in Table 4
for $c=10$. The values are the same for $c=5$ up to $t=14$ after
which they decrease more slowly. We compare $m$ with
$m_\infty=(1+4t/b)^{-1/2}$. This is the `same' constant calculated
using Fourier transforms when $a=\infty$, i.e. for the semigroup
on $L^2(\R)$ when the initial function is $f(x)=\rme^{-x^2}$. The
two agree up to $t=10$, which is all that one could expect. All of
the results confirm that the pseudospectral approximation to the
semigroup is highly reliable for the stated values of $a$ and $b$,
at least for this choice of the initial function $f$.

\begin{table}[h]  
\[
\begin{array}{ccc}
t &  m &m_\infty\\
0&  1.0000 &1.0000 \\
2&  0.8451 &0.8452\\

4&  0.7454  &  0.7454 \\

6&  0.6742 &  0.6742\\

8&  0.6202 &  0.6202\\

10&   0.5593 & 0.5774\\

12&  0.1268 & 0.5423\\

14&   0.0049 &0.5130\\

16&   0.0000 & 0.4880\\
\end{array}
\]
\begin{center}
Table 4
\end{center}
\end{table}

We repeated the calculations leading to Figure 1, but with the
initial function $g(x)=1$ for all $x\in [0,a]$. This is a much
more serious test of the method since $g$ does not satisfy the
boundary conditions even approximately. With $N=15$ and $c=5$ we
obtained the results shown in Figure 2. One sees that $g_t$ is
close to the characteristic function of $[0,a-t]$, but smoothed
out because of the diffusion term in $A$. By contrast with the
similar calculation in Section 4, there is no Gibbs phenomenon,
presumably again because of the diffusion term. For smaller values
of $b$, such as $b=5$, the fact that $g_t(0)=0$ for all $t>0$ is
much more obvious.

\begin{figure}[h]  
\begin{picture}(200,350)(0,0)
\includegraphics{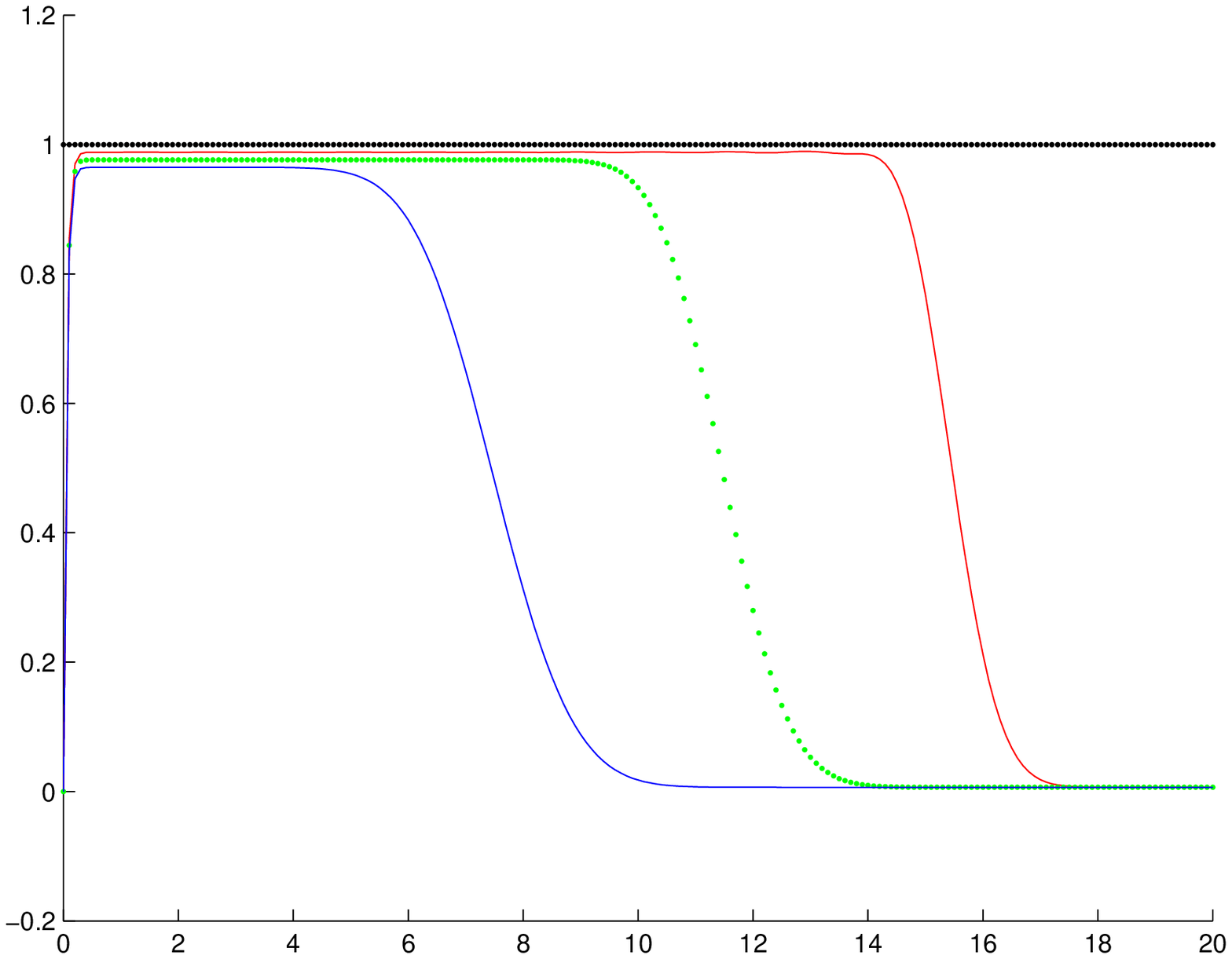}
\end{picture}
\begin{center}
Figure 2. Graphs of $g$, $g_8$ (dotted) and $g_4$, $g_{12}$
(solid)
\end{center}
\end{figure}
\clearpage

Table 5 lists the first few eigenvalues $\lam_n$ and approximate
eigenvalues $\mu_s$ of $A$ in decreasing order of their real
parts, where $a=20$, $b=20$ and $c=5$. The largest eigenvalue
$-5.001$ controls the asymptotic decay of the semigroup as
$t\to\infty$, but it has little influence on the size of $\norm
T_tf \norm$ for $t=10$. One of the main reasons for the accuracy
of the pseudospectral expansion is the fact that there are so many
approximate eigenvalues whose real parts are close to zero. For
$c=10$ the real parts of these $\mu_s$ decrease from $-0.488$ to
$-0.729$.

\begin{table}[h]  
\[
\begin{array}{cc}
\lam_n & \mu_s \\
-5.001&-0.247   \\
-5.005& -0.252\pm 0.306i \\

-5.011& -0.252\pm 0.306i  \\

-5.020& -0.267\pm 0.613i \\

-5.031& -0.291\pm 0.919i \\

-5.044& -0.326\pm 1.225i  \\

-5.060& -0.370\pm 1.532i \\

-5.079& -0.425\pm 1.838i \\

\end{array}
\]
\begin{center}
Table 5
\end{center}
\end{table}

\vspace{2 cm}

{\bf Acknowledgements} We would like to acknowledge financial
support under the EPSRC grant GR/R81756/01.

\vspace{2 cm}

Department of Mathematics \\
King's College\\
Strand\\
London\\
WC2R 2LS\\
England

E.Brian.Davies@kcl.ac.uk

\end{document}